\documentclass[12pt]{amsart}
\usepackage{amssymb,latexsym,amsmath,enumerate,nicefrac}
\usepackage{amsaddr}

\usepackage[total={6.5in,8.5in}, top=1.25in, left=1in]{geometry}

\usepackage{color}
\definecolor{myurlcolor}{rgb}{0.6,0,0}
\definecolor{mycitecolor}{rgb}{0,0,0.8}
\definecolor{myrefcolor}{rgb}{0,0,0.8}
\usepackage[bookmarks=false,pagebackref=true]{hyperref}
\hypersetup{colorlinks,
linkcolor=myrefcolor,
citecolor=mycitecolor,
urlcolor=myurlcolor}

\newcommand{\vertiii}[1]{{\vert\kern-0.25ex\vert\kern-0.25ex\vert #1 
    \vert\kern-0.25ex\vert\kern-0.25ex\vert}}

\newcommand{\bvertiii}[1]{{\big\vert\kern-0.25ex\big\vert\kern-0.25ex\big\vert #1 
    \big\vert\kern-0.25ex\big\vert\kern-0.25ex\big\vert}}

\newcommand{\Bvertiii}[1]{{\Big\vert\kern-0.25ex\Big\vert\kern-0.25ex\Big\vert #1 
    \Big\vert\kern-0.25ex\Big\vert\kern-0.25ex\Big\vert}}    
    
\newtheorem{theorem}{Theorem}[section]
\newtheorem{lemma}[theorem]{Lemma}

\newtheorem{proposition}[theorem]{Proposition}
\newtheorem*{corollary}{Corollary}
\theoremstyle{definition}

\theoremstyle{remark}

\usepackage[all]{xy}
\usepackage{comment}
 \usepackage{tikz}
\usetikzlibrary{arrows}
 \usetikzlibrary{matrix,positioning}
\usetikzlibrary{calc}
\usepackage{graphics}
\usepackage{dsfont}
\usepackage{amsmath}
\usepackage{amsfonts,amssymb}
\usepackage{mathrsfs}
\usepackage{mathtools}
\usepackage{comment}

\allowdisplaybreaks

\makeatletter
\newcommand*\bigcdot{\mathpalette\bigcdot@{.6}}
\newcommand*\bigcdot@[2]{\mathbin{\vcenter{\hbox{\scalebox{#2}{$\m@th#1\bullet$}}}}}
\makeatother

\newcommand{\Q}{\mathbb{Q}}

\newcommand{\e}{{\rm e}}

\renewcommand{\(}{\left(} 
\renewcommand{\)}{\right)}

\usepackage{stmaryrd}

\pgfkeys{/tikz/.cd, arrowhead/.default=-1pt, arrowhead/.code={}}

\begin{document}

\title{Components and Exit Times of Brownian Motion in two or more $p$-Adic Dimensions}

%

\author{Rahul Rajkumar$^1$ \and David Weisbart$^2$}
\address{
\begin{tabular}[h]{cc}
 $^{1,2}$Department of Mathematics\\
  University of California, Riverside
  \end{tabular}
  }
\email{$^1$rajkumar@math.ucr.edu} \email{$^2$weisbart@math.ucr.edu}



\maketitle

\pagestyle{plain}


\begin{abstract}
The fundamental solution of a pseudo-differential equation for functions defined on the $d$-fold product of the $p$-adic numbers, $\mathds{Q}_p$, induces an analogue of the Wiener process in $\mathds{Q}_p^d$. As in the real setting, the components are $1$-dimensional $p$-adic Brownian motions with the same diffusion constant and exponent as the original process. Asymptotic analysis of the conditional probabilities shows that the vector components are dependent for all time. Exit time probabilities for the higher dimensional processes reveal a concrete effect of the component dependency.
\end{abstract}


\tableofcontents

\section{Introduction}\label{sec:intro}

Two main ideas initially motivated the study of non-Archimedean diffusion: the idea that non-Archimedean physical models describe the observed ultrametricity in certain complex systems, and the idea that the extremely small scale structure of spacetime could be non-Archimedean.  Ultrametric structures in spin glasses were already implicit in Parisi's investigations in \cite{Parisi:PRL:1979} and \cite{Parisi:JPA:1980}.  In \cite{Volovich:1987}, Volovich proposed the idea that spacetime could have a non-Archimedean structure at small enough distance and time scales and he initiated a program to study analogues of physical theories with $p$-adic state spaces for this reason. Varadarajan discussed this contribution of Volovich in \cite[Chapter 6]{VSV:Reflections:2011}. The study of diffusion in non-Archimedean local fields goes back more than 30 years, and the study of diffusion processes in vector spaces over such fields goes back nearly as far.  Seminal works in this area include \cite{koch92} and \cite{alb}.  In \cite{koch92}, Kochubei gave the fundamental solution to the $p$-adic analogue of the diffusion equation, developed a theory of $p$-adic diffusion equations, and proved a Feynman-Kac formula for the operator semigroup with a $p$-adic Schr\"{o}dinger operator as its infinitesimal generator.  Albeverio and Karwowski constructed in \cite{alb} a continuous time random walk on $\mathds Q_p$, computed its transition semigroup and infinitesimal generator, and showed among other things that the associated Dirichlet form is of jump type.

For any finite dimensional vector space $\mathcal S$ with coefficients in a division algebra that is finite dimensional over a non-Archimedean local field of arbitrary characteristic, Varadarajan constructed in \cite{var97} a general class of diffusion processes with sample paths in the Skorohod space $D([0, \infty)\colon\mathcal S)$ of c\`adl\`ag paths that take values in $\mathcal S$. The current work follows the approach of \cite{var97} and takes it as a starting point, but specializes to the setting where, for any prime number $p$ and any natural number $d$, paths take values in the $d$-fold Cartesian product of the $p$-adic numbers, $\mathds Q_p^d$.  The results of \cite{var97} show that there is a triple $\big(D([0, \infty)\colon\mathds Q_p^d), P^d, \vec{X}\big)$ so that $P^d$ is a probability measure on $D([0, \infty)\colon\mathds Q_p^d)$ and, for any $\omega$ in $D([0, \infty):\mathds Q_p^d)$ and any positive $t$, \[\vec{X}(t, \omega) = \omega(t).\] Furthermore, the probability measure on $\mathds Q_p^d$ that for any Borel set $B$ of $\mathds Q_p^d$ is given by \[B\mapsto P^d(\vec{X}(t, \omega) \in B)\] has a density function that is a solution to a pseudo-differential equation that is analogous to the real diffusion equation.  Refer to any process of this type as a \emph{Brownian motion in $\mathds Q_p^d$}.

For any Brownian motion $\vec{X}$ in $\mathds Q_p^d$ with sample paths in $(D([0, \infty)\colon\mathds Q_p^d), P^d)$, the current paper establishes that the component processes of $\vec{X}$ are each Brownian motions in $\mathds Q_p$ with the same parameters (diffusion constant and exponent) as $\vec{X}$, and that for no positive real $t$ are the components of $\vec{X}_t$ independent.  Section~\ref{sec:NormProd} briefly reviews some necessary results in $p$-adic analysis, discusses the max-norm process that is a special case of the more general process that Varadarajan discusses in \cite{var97}, and determines the infinitesimal generator of the \emph{product process} in $\mathds Q_p^d$.  Section~\ref{sec:components} studies the component processes of the max-norm process.  The main results are Theorems~\ref{thm:components:marginal_distributions} and \ref{theorem:components:epsilon}.  Theorem~\ref{thm:components:marginal_distributions} establishes that the component processes are Brownian motions in $\mathds Q_p$.  Theorem~\ref{theorem:components:epsilon} gives precise estimates on certain conditional probabilities that establish the dependency of the component processes.  The effect of the dependency of the components becomes strikingly apparent in the calculation in Section~\ref{sec:exit} of the exit probabilities for the product and max-norm process that generalize \cite[Theorem 3.1]{Weisbart:2021}.

Dragovich, Khrennikov, Kozyrev, Volovich, and Zelenov give a detailed review of the history of the research in non-Archimedean mathematical physics in \cite{DKKVZ:2017} that updates the earlier review \cite{DKKV:2009} by the first four authors.  This review helps to put the current paper in context and discusses many areas where the current paper could find application.  In their recent book \cite{KKZ:2018}, Khrennikov, Kozyrev, and Z\'{u}\~{n}iga-Galindo discuss many applications of ultrametric pseudodifferential equations, including many interesting recent developments.  This work and the references therein also present many areas where the current paper could be useful.  The works \cite{ave3, ave4, ave} of Avetisov, Bikulov, Kozyrev, and Osipov that deal with $p$-adic models for complex systems seem to be particularly relevant to this current paper, as is the work \cite{Avetisov-Bikulov:PSIM:2009} of Avetisov and Bikulov that involve biological applications. Ultrametricity can be found in data structures, as Bradley discusses in \cite{Bradley:2017}.  The current paper may find application in the study of data structures and, in particular, in the recent work \cite{Bradley-Keller-Weinmann: 2018} of Bradley, Keller, Weinmann, as well as in Bradley's work, \cite{Bradley: 2019}, and in the work \cite{Bradley-Jahn: 2022}, of Bradley and Jahn.



\section{The Norm and Product Processes}\label{sec:NormProd}

See Gouv\^{e}a's book \cite{Gov} for an accessible supplement to the cursory review of $p$-adic analysis that the current section presents.  For more detail, see the book \cite{Ram} of Ramakrishnan and Valenza, and Weil's book \cite{Weil}.  Vladimirov, Volovich, and Zelenov give a self-contained introduction to $p$-adic analysis and mathematical physics in their now classic book \cite{vvz}.  Z\'{u}\~{n}iga-Galindo's recent book \cite{Zuniga:PDiffEBook} is an accessible and current reference that, among other things, investigates diffusion processes of a type that includes the max-norm processes that appear shortly. 

This section follows the presentation of \cite{Weisbart:2021}, but generalizes it to the higher dimensional $p$-adic setting that is necessary for the sections that follow.

\subsection{Analysis in $\mathds Q_p^d$}\label{sec:NormProd:sub:analysis}

For any prime $p$, denote by $|\cdot|$ the absolute value on $\mathds Q_p$.  For any natural number $d$ and any $d$-tuple $(x_1, \dots, x_d)$ in $\mathds Q_p^d$, write \[\vec{x} = (x_1, \dots, x_d).\]  Denote by $\|\cdot\|$ the \emph{max-norm} on $\mathds Q_p$ that takes any $\vec{x}$ in ${\mathds Q}_p^d$ to $\|\vec{x}\|$, where \[\|\vec{x}\| = \max_{i\in\{1, \dots, d\}}|x_i|.\] The max-norm induces an ultrametric on $\mathds Q_p^d$.  The general linear group in $d$ dimensions with coefficients in $\mathds Z\slash p\mathds Z$, $GL_d(\mathds Z\slash p\mathds Z)$, is the maximal compact subgroup of the $p$-adic general linear group in $d$ dimensions, $GL_d(\mathds Q_p)$.  Since the real orthogonal group in $d$ dimensions, $O_d(\mathds R)$, is the maximal compact subgroup of the real general linear group $GL_d(\mathds R)$, the group $GL_d(\mathds Z\slash p\mathds Z)$ is the natural analogue in the $\mathds Q_p^d$ setting of the real orthogonal group $O_d(\mathds R)$.  The max-norm is $GL_d(\mathds Z\slash p\mathds Z)$ invariant, and so the max-norm is a natural analogue in the $\mathds Q_p^d$ setting of the Euclidean norm on $\mathds R^d$.

For each $\vec{x}$ in $\mathds Q_p^d$, denote respectively by $B_d(k,\vec{x})$ and $S_d(k,\vec{x})$ the ball and the circle of radius $p^k$ with center at $\vec{x}$, the compact open sets \[B_d(k, \vec{x}) = \big\{\vec{y}\in \mathds Q_{p}^d\colon \|\vec{y}-\vec{x}\| \leq p^k\big\} \quad {\rm and}\quad {S}_d(k, \vec{x}) = \big\{\vec{y}\in \mathds Q_{p}^d\colon \|\vec{y}-\vec{x}\| = p^k\big\}.\]  The unit ball in $\mathds Q_{p}^d$ is the $d$-fold Cartesian product, $\mathds Z_{p}^d$, of the ring of integers in $\mathds Q_p$.  For balls and circles with centers at the origin, simplify notation by respectively denoting by $B_d(k)$ and $S_d(k)$ the ball and the circle of radius $p^k$ in $\mathds Q_p^d$ with center $\vec{0}$.  To further simplify notation, suppress $d$ in the notation for balls and circles to indicate that $d$ is equal to 1.

For any $x$ in $\mathds Q_{p}$, there is a unique function \begin{align}\label{theafunction}a_{x}\colon \mathds Z\to \{0,1,\dots, p-1\}\end{align} with the property that \[x = \sum_{k\in\mathds Z}a_x(k)p^k.\]   Denote by $\{x\}$ the \emph{fractional part of $x$}, the sum %
\begin{equation}\label{eq:NormProd:def:fracpart}
\{x\} = \sum_{k<0} a_{x}(k)p^k.
\end{equation}%
For any $x$ in $\mathds Q_p$, the support of $a_x$ is bounded below and so the sum that defines $\{x\}$ in \eqref{eq:NormProd:def:fracpart} is a finite sum. Take $\chi$ to be the additive character on $\mathds Q_p$ that is given by \[\chi(x) = {\rm e}^{2\pi{\sqrt{-1}}\{x\}}.\] For any natural number $d$, $\mathds Q_p^d$ is a totally disconnected, self-dual, locally compact, Hausdorff abelian group.  The additive character $\chi$ induces an isomorphism between $\mathds Q_p^d$ and its Pontryagin dual, $\big(\mathds Q_p^d\big)^\ast$.  Namely, for any additive character $\phi$ in $\big(\mathds Q_p^d\big)^\ast$, there is a $\vec{y}$ in $\mathds Q_p^d$ so that for any $\vec{x}$ in $\mathds Q_p^d$, \[\phi(\vec{x}) = \chi(\vec{x}\cdot \vec{y}) \quad \text{where}\quad \vec{x}\cdot \vec{y} = x_1y_1 +\cdots + x_dy_d.\]  Take $\mu_d$ to be the unique Haar measure on $\mathds Q_{p}^d$ for which $\mathds Z_{p}^d$ has unit measure.  For any integer $k$, the substitution formula for $p$-adic integrals implies that $\mu_1(B(k))$ is equal to $p^k$.  The measure $\mu_d$ is a product measure and $B_d(k)$ is the $d$-fold Cartesian product of the ball $B(k)$. Furthermore, the circle $S(k)$ is the set $B(k)\setminus B(k-1)$.  Translation invariance of $\mu_d$ therefore implies that for any $\vec{x}$ in $\mathds Q_p^d$,
\begin{equation}\label{EQ:NormProd:Ball_Sphere_measure_d}
\mu_d(B_d(k, \vec{x})) = p^{kd} \quad {\rm and}\quad \mu_d(S_d(k, \vec{x})) = p^{kd}- p^{(k-1)d} = p^{kd}\left(1 - \tfrac{1}{p^d}\right).
\end{equation}

Initially define $\mathcal F_d$ and $\mathcal F_d^{-1}$ on $L^1(\mathds Q_{p}^d)$ by \[(\mathcal F_df)(\vec{y}) = \int_{\mathds Q_{p}^d}\chi(-\vec{x}\cdot \vec{y})f(\vec{x})\,{\rm d}\mu_d\!\(\vec{x}\) \quad {\rm and}\quad (\mathcal F_d^{-1}f)(\vec{y}) = \int_{\mathds Q_{p}^d}\chi(\vec{x}\cdot \vec{y})f(\vec{x})\,{\rm d}\mu_d\!\(\vec{x}\).\]  These operators are unitary on $L^1(\mathds Q_{p}^d)\cap L^2(\mathds Q_{p}^d)$ and extend to unitary operators on $L^2(\mathds Q_{p}^d)$.  The extensions of these operators, again denoted by $\mathcal F_d$ and $\mathcal F_d^{-1}$, are the Fourier and inverse Fourier transforms on $L^2(\mathds Q_{p}^d)$, respectively.  To simplify notation, suppress the measure in the notation for integrals by writing ${\rm d}\vec{x}$ to mean ${\rm d}\mu_d\big(\vec{x}\,\big)$.  Furthermore, write ${\rm d}x$ to mean ${\rm d}\mu_1(x)$ in the case when $d$ is equal to $1$ and the integral is taken over a subset of $\mathds Q_p$.  

The following lemma is helpful for performing calculations that involve the integration of characters.

\begin{lemma}\label{lem:NormProd:CharInt}
For any $m$ and $n$ in ${\mathds Z}$, 
\begin{equation*}
\int_{B_d(m)} \int_{B_d(n)} \chi(\vec{x}\cdot \vec{y}\,) \,{\rm d}\vec{x} \,{\rm d}\vec{y} = p^{d(n + \min(-n,m))}.
\end{equation*}
\end{lemma}

\begin{proof}
Since the character $\chi$ is identically equal to 1 on $\mathds Z_p$ and the sum of the $p^{\rm th}$ roots of unity is equal to $0$, 
\begin{equation}\label{lemma:NormProd:BasicCharInt}
\int_{B(n)} \chi(x)\,{\rm d}x = \begin{cases}p^{n}&\mbox{if }p^{n}\leq 1\\0&\mbox{if }p^{n}> 1.\end{cases}
\end{equation}
For any $y$ in $\mathds Q_p$, the substitution formula for $p$-adic integration implies that 
\begin{equation}\label{lemma:NormProd:BasicCharInt}
\int_{B(n)} \chi(xy)\,{\rm d}x = \begin{cases}p^{n}&\mbox{if }p^{n}\leq \frac{1}{|y|}\\0&\mbox{if }p^{n}> \frac{1}{|y|}.\end{cases}
\end{equation}
For any subset $S$ of $\mathds Q_p^d$, take ${\mathds 1}_S$ to be the indicator function on $S$.   Equation~\ref{lemma:NormProd:BasicCharInt} implies that 
\begin{equation}\label{lemma:NormProd:BasicCharInt:Indicator}
\int_{B_d(n)} \chi(\vec{x}\cdot\vec{y}\,)\,{\rm d}\vec{x} = p^{dn}{\mathds 1}_{B_d(-n)}(\vec{y}\,),
\end{equation}
and so
\begin{align*}
\int_{B_d(m)} \int_{B_d(n)} \chi(\vec{x}\cdot\vec{y}\,) \,{\rm d}\vec{x}\, {\rm d}\vec{y} &= \int_{B_d(-m)} p^{dn}{\mathds 1}_{B_{-n}}(\vec{y}\,) \,{\rm d}\vec{y} \\&= p^{d(n + \min(-n,m))}.
\end{align*}
\end{proof}

\subsection{The Max-Norm Process}\label{sec:NormProd:sub:1-d}

Denote by $SB(\mathds Q_{p}^d)$ the \emph{Schwartz-Bruhat} space of complex valued, compactly supported, locally constant functions on $\mathds Q_{p}^d$.  This space is the $\mathds Q_p^d$ analogue of the space of complex valued, compactly supported, smooth functions on $\mathds R^d$.  Unlike its real analogue, $SB(\mathds Q_{p}^d)$ is closed under the Fourier transform.  For any positive real number $b$, take ${\mathcal M}_b$ to be the multiplication operator that acts on $SB(\mathds Q_{p}^d)$ by \[({\mathcal M}_bf)(\vec{x}) = \|\vec{x}\|^bf(\vec{x}).\]  Take $\Delta_{b,d}$ to be the self-adjoint closure of the densely defined operator that acts on $SB(\mathds Q_{p}^d)$ by \begin{equation}\label{EQ:NormProd:pseudoDelta}\big(\Delta_{b,d} f\big)(\vec{x}) = \big(\mathcal F^{-1}_d{\mathcal M_b}\mathcal F_df\big)\!(\vec{x}).\end{equation}  Extend $\Delta_{b,d}$ to act on a complex valued function $f$ on $\mathds R_+\times \mathds Q_{p}^d$ by currying variables, as in \cite{Weisbart:2021}.  This is to say that for each $t$ in $\mathds R_+$, if the function $f(t, \cdot)$ is in the domain of $\Delta_{b,d}$, then $f$ is in the domain of the extended operator, again denoted by $\Delta_{b,d}$, and that \[(\Delta_{b,d}f)(t,\vec{x}) := (\Delta_{b,d}f(t, \cdot))(\vec{x}).\]  This extension is the \emph{Taibleson-Vladimirov operator with exponent} $b$. To simplify notation, suppress $d$ to denote that $d$ is equal to $1$. Again follow \cite{Weisbart:2021} by currying variables to extend the Fourier and inverse Fourier transforms to act on functions on $\mathds R_+\times \mathds Q_{p}^d$.  For any positive real number $\sigma$, the pseudo-differential equation \begin{align}\label{EQ:NormProd:DiffusionEQ} \dfrac{{\rm d}f(t,\vec{x})}{{\rm d}t} = -\sigma\Delta_{b,d} f(t,\vec{x})\end{align} is a \emph{$d$-dimensional $p$-adic diffusion equation} and has fundamental solution $\rho_d$, where for each $t$ in $(0,\infty)$ and for each $\vec{x}$ in $\mathds Q_p^d$, \[\rho_d(t,\vec{x}) = \left(\mathcal F_d^{-1}\e^{-\sigma t\|\cdot\|^b}\right)\!(\vec{x}).\] Once again, suppress $d$ in the notation for $\rho_d$ when $d$ is equal to $1$.  %

With the necessary modifications to include the diffusion constant $\sigma$, follow the arguments in \cite{var97} to see that $\rho_d(t,\cdot)$ is a probability density function that gives rise to a probability measure $P^d$ on $D([0,\infty) \colon \mathds Q_{p}^d)$ that is concentrated on the set of paths originating at 0.  The inclusion of a diffusion constant that is not equal to 1 amounts to a rescaling of the real time parameter.   The \emph{max-norm process} $\vec{X}$ is the stochastic process that, for any pair $(t, \omega)$ in $[0,\infty)\times D([0,\infty) \colon \mathds Q_{p}^d)$, is defined by \[\vec{X}(t, \omega) = \omega(t).\]   This stochastic process specializes the process that Varadarajan constructed in \cite{var97} to the setting of $\mathds Q_p^d$.   If $d$ is equal to $1$, write $X$ rather than $\vec{X}$.  The process $X$ is the process discussed in \cite{Weisbart:2021} and is the \emph{Brownian motion} in $\mathds Q_p$ with \emph{diffusion constant} equal to $\sigma$ and \emph{diffusion exponent} equal to $b$.

Denote by $\vec{X}_t$ the random variable $\vec{X}(t, \cdot)$ that takes any path $\omega$ in $D([0,\infty) \colon \mathds Q_{p}^d)$ to the value $\omega(t)$.  The density function $\rho_d(t,\cdot)$ for $\vec{X}_t$ satisfies the equality
\begin{align}\label{pdf:NormProcess:d}
\rho_d(t,\vec{x}) & = \int_{\mathds Q_p^d} \chi(\vec{x}\cdot \vec{y})\e^{-\sigma t\|\vec{x}\|^b}\,{\rm d}\vec{y}\notag\\&= \sum_{r\in \mathds Z} \e^{-\sigma tp^{rb}}\int_{S_d(r)} \chi(\vec{x}\cdot \vec{y})\,{\rm d}\vec{y}\notag\\ &= \sum_{r\in \mathds Z} \Big(\e^{-\sigma tp^{rb}} - \e^{-\sigma tp^{(r+1)b}}\Big)\int_{B_d(r)} \chi(\vec{x}\cdot \vec{y})\,{\rm d}\vec{y}\notag\\& = \sum_{r\in\mathds Z}\Big(\e^{-\sigma t p^{rb}} - \e^{-\sigma  t p^{(r+1)b}}\Big)p^{dr}{\mathds 1}_{B_d(-r)}(\vec{x}).
\end{align}

\subsection{The Product Process}\label{sec:NormProd:higher-d}

A function $f$ is a \emph{Schwartz-Bruhat monomial} with domain ${\mathds Q}_p^d$ if for every $i$ in $\{1, \dots, d\}$ there is a Schwartz-Bruhat function $f_i$ on $\mathds Q_p$ so that for every $\vec{x}$ in ${\mathds Q}_p^d$,
\begin{equation} 
\vec{x} = (x_1, \dots, x_d) \quad \text{implies that}\quad f(\vec{x}) = \prod_{i=1}^d f_i(x_i).
\end{equation}
The space $SB_0(\mathds Q_{p}^d)$ of \emph{simple Schwartz-Bruhat} functions on $\mathds Q_p^d$ is the $d$-fold algebraic tensor product of the space of functions $SB({\mathds Q}_p)$---It is the complex vector space of functions that are finite linear combinations of Schwartz-Bruhat monomials.  Each simple Schwartz-Bruhat function is a function in $L^2({\mathds Q}_p^d)$, and $L^2({\mathds Q}_p^d)$ is the analytic completion of $SB_0(\mathds Q_{p}^d)$ under the $L^2$-norm on $L^2({\mathds Q}_p^d)$.

For each $i$ in $\{1, \dots, d\}$, take $H_i$ to be the linear extension to $SB_0(\mathds Q_{p}^d)$ of the operator that is initially defined for any Schwartz-Bruhat monomial $f$ by \[H_i(f) = f_1 \otimes \cdots \otimes f_{i-1} \otimes \sigma\Delta_b f_i\otimes f_{i+1} \otimes \cdots \otimes f_d.\] Take $H_0$ to be the sum \[H_0 = H_1 + \cdots + H_d.\]  The operator $H_0$ is the Fourier transform of real valued multiplication operator on a dense subset of $L^2({\mathds Q}_p^d)$, and so it is essentially self adjoint on $SB_0(\mathds Q_{p}^d)$.  Take $H$ to be the self adjoint closure of $H_0$ on $L^2({\mathds Q}_p^d)$.

The max-norm process $\vec{X}$ has a radially symmetric law with respect to the max-norm, but it is not the only possible choice for a Brownian motion in ${\mathds Q}_p^d$.  Take $P$ to be the probability measure on $D([0,\infty) \colon {\mathds Q}_p)$ such that the triple $(D([0,\infty) \colon {\mathds Q}_p), P, X)$ is the 1-dimensional Brownian motion studied in \cite{Weisbart:2021} with the property that, for each positive $t$, the probability density function $\rho(t,\cdot)$ for $X_t$ satisfies \eqref{EQ:NormProd:DiffusionEQ} with $d$ equal to $1$.  For each positive real number $t$, take $\vec{Y}_t$ to be the random variable that, for any path $\omega$ in the $d$-fold Cartesian product $D([0,\infty) \colon \mathds Q_{p})^d$, is given by \[\vec{Y}_t(\omega) = \omega(t) = (\omega_1(t), \dots, \omega_d(t)).\]  Take $\vec{Y}$ to be the function that is defined for all $t$ in $(0,\infty)$ by \[\vec{Y}(t) = \vec{Y}_t.\]  The process that is given by the triple $\big(D([0,\infty) \colon {\mathds Q}_p^d), \otimes_{i=1}^d P, \vec{Y}\big)$ is the \emph{product process}, which is radially symmetric if and only if $d$ is equal to $1$.  Simplify notation by henceforth writing $\otimes^dP$ rather than $\otimes_{i=1}^dP$.

For any positive real $t$, take $g_d(t, \cdot)$ to be the probability density function for $\vec{Y}_t$.  The independence of the components of $\vec{Y}$ together with the product rule for differentiation implies that, for any $\vec{x}$ in ${\mathds Q}_p^d$, 
\begin{align*}
\frac{{\rm d}}{{\rm d}t}g_d(t,\vec{x}) & = \frac{{\rm d}}{{\rm d}t}\big(\rho(t,x_1)\rho(t,x_2)\cdots \rho(t,x_d)\big)\\ & = \sigma H g_d(t, \vec{x}),
\end{align*}
and so $H$ is the infinitesimal generator of the product process on ${\mathds Q}_p$.  

In the setting of Brownian motion in ${\mathds R}^d$, where the norm is the usual Euclidean $\ell^2$-norm, the analogous norm process is that which is given by the product of $d$ independent Brownian motions in ${\mathds R}$ that have the same diffusion constant.  The next section will show that this is not the case in the ${\mathds Q}_p$ setting, and the goal of the present work is to precisely understand some differences between the stochastic processes $\vec{X}$ and $\vec{Y}$ when $d$ is greater than 1.


\section{The Component Processes and their Dependence}\label{sec:components}

Section~\ref{sec:components:sub:marginals} establishes that the processes $\vec{X}$ and $\vec{Y}$ are similar in that the components of each are themselves Brownian motions with the same diffusion constant, $\sigma$, and the same diffusion exponent, $b$.  Section~\ref{sec:components:sub:dependency} establishes that the max-norm process, $\vec{X}$, has dependent components, where the product process, $\vec{Y}$, has independent components, and so $\vec{X}$ and $\vec{Y}$ are qualitatively different.

\subsection{Calculation of the Marginals}\label{sec:components:sub:marginals}

For any natural number $i$ in $\{1, \dots, d\}$, take $X^{(i)}$ to be the stochastic process that is given by the $i^{\rm th}$ component of the max-norm process $\vec{X}$.  The goal of this subsection is to prove Theorem~\ref{thm:components:marginal_distributions}.

Some additional notation facilitates the presentation.  For each $\vec{x}$ in $\mathds Q_p^d$, denote by $\vec{x}_i$ the vector in $\mathds Q_p^{d-1}$ that is given by \[\vec{x}_i = (x_1, \dots, x_{i-1}, x_{i+1}, \dots, x_d).\]  To simplify notation, identify any ordered pair $(x_1, (x_2, \dots, x_d))$ in $\mathds Q_p\times \mathds Q_p^{d-1}$ with the $d$-tuple $(x_1, \dots, x_d)$ in $\mathds Q_p^d$.  For each positive real number $t$, take $\rho^{(i)}(t, \cdot)$ to be the probability density function for the $i^{\rm th}$ component function of $\vec{X}_t$, the random variable $X_t^{(i)}$.  For each $x$ in $\mathds Q_p$, the marginal $\rho^{(i)}(t, x)$ is given by the integral  %
\begin{align}\label{EQ:components_marginals:motivation}
\rho^{(i)}(t, x) & = \int_{{\mathds Q}_p^{d-1}} \rho_d(t, (x_1, \dots, x_{i-1}, x, x_{i+1}, \dots, x_d)) \,{\rm d}\vec{x}_i \notag\\ 
& = \int_{{\mathds Q}_p^{d-1}} \rho_d(t, (x, \vec{x}_1)) \,{\rm d}\vec{x}_1\notag\\ 
&=  \int_{{\mathds Q}_p^{d-1}}  \int_{{\mathds Q}_p^d} \chi((x, \vec{x}_1)\cdot\vec{y}\,)\e^{-\sigma t\|\vec{y}\|^b} \,{\rm d}\vec{y} \,{\rm d}\vec{x}_i,
\end{align}
where the change of variables formula and the symmetry of the integrand together imply the penultimate equality.

Switching the order of integration in \eqref{EQ:components_marginals:motivation} that determines the marginals facilitates calculation of the law for the component processes. The failure of the integrand to be absolutely integrable with respect to the product measure on ${\mathds Q}_p^{d-1} \times {\mathds Q}_p^d$ precludes a naive application of the Fubini-Tonelli theorem that would quickly verify Theorem~\ref{thm:components:marginal_distributions}.  The technical complication necessitates the following more involved argument and implicitly makes use of the fact that the given integral is an oscillatory integral.

\begin{theorem}\label{thm:components:marginal_distributions}
Each $X^{(i)}$ is a Brownian motion in $\mathds Q_p$ with diffusion constant equal to $\sigma$ and diffusion exponent equal to $b$.
\end{theorem}

\begin{proof}
For any $i$ in $\{1, \dots, d\}$, the radial symmetry of $\rho_d$ and \eqref{EQ:components_marginals:motivation} together imply that $X^{(i)}$ and $X^{(1)}$ have the same law.  For any positive $t$, the characteristic function $\phi_d(t, \cdot)$ of $\vec{X}_t$ is given for each $\vec{y}$ in $\mathds Q_p^d$ by \[\phi_d(t, \vec{y}\,) = \e^{-\sigma t\|\vec{y}\,\|^b}.\]  The function $\phi_d(t, \cdot)$ is bounded and integrable, and so the Fubini-Tonelli theorem guarantees that
\begin{align*}
\rho^{(1)}(t, x_1) & = \int_{{\mathds Q}_p^{d-1}} \rho_d(t, (x_1, \vec{x}_1))\,{\rm d}\vec{x}_1\\
& = \int_{{\mathds Q}_p^{d-1}}\left\{\int_{{\mathds Q}_p^{d}} \chi((x_1, \vec{x}_1)\cdot \vec{y}\,)\e^{-\sigma t\|\vec{y}\,\|^b}\,{\rm d}\vec{y}\right\}{\rm d}\vec{x}_1\\%
& = \int_{{\mathds Q}_p^{d-1}}\left\{\int_{{\mathds Q}_p} \int_{{\mathds Q}_p^{d-1}} \chi(x_1y_1)\chi(\vec{x}_1\cdot \vec{y}_1)\e^{-\sigma t\|(y_1, \vec{y}_1)\|^b}\,{\rm d}\vec{y}_1{\rm d}y_1\right\}{\rm d}\vec{x}_1.
\end{align*}
For each $y_1$, decompose the innermost integral into a sum of integrals over the ball of radius $|y_1|$ and its complement to obtain the equalities%
\begin{align*}
\rho^{(1)}(t, x_1) & = \int_{{\mathds Q}_p^{d-1}}\left\{\int_{\mathds Q_p}\int_{\|\vec{y}_1\| \leq |y_1|} \chi(x_1y_1)\chi(\vec{x}_1\cdot \vec{y}_1)\e^{-\sigma t\|(y_1, \vec{y}_1)\|^b}\,{\rm d}\vec{y}_1\,{\rm d}y_1\right\}{\rm d}\vec{x}_1 \\&\qquad\qquad\qquad+ \int_{{\mathds Q}_p^{d-1}}\left\{\int_{\mathds Q_p}\int_{\|\vec{y}_1\|>|y_1|} \chi(x_1y_1)\chi(\vec{x}_1\cdot \vec{y}_1)\e^{-\sigma t\|(y_1, \vec{y}_1)\|^b}\,{\rm d}\vec{y}_1\,{\rm d}y_1\right\}{\rm d}\vec{x}_1\\
& = \int_{{\mathds Q}_p^{d-1}}\Bigg\{\int_{\mathds Q_p}\chi(x_1y_1)\Bigg[\e^{-\sigma t|y_1|^b}\int_{\|\vec{y}_1\| \leq |y_1|} \chi(\vec{x}_1\cdot \vec{y}_1)\,{\rm d}\vec{y}_1 \\&\qquad\qquad\qquad+ \int_{\|\vec{y}_1\|>|y_1|} \chi(\vec{x}_1\cdot \vec{y}_1)\e^{-\sigma t\|\vec{y}_1\|^b}\,{\rm d}\vec{y}_1\Bigg]\,{\rm d}y_1\Bigg\}{\rm d}\vec{x}_1.
\end{align*}

Since $\rho_d(t, \cdot)$ is integrable, %
\begin{equation}\label{EQ:components_marginals:Rhod_as_a_Limit}\rho^{(1)}(t, x_1) = \lim_{n\to\infty}\int_{B_{d-1}(n)} \rho_d(t, (x_1, \vec{x}_1))\,{\rm d}\vec{x}_1.\end{equation}%
Take $I_1(x_1)$ and $I_2(x_1)$ to be the quantities %
\begin{equation}\label{EQ:components_marginals:I1a}I_1(x_1) = \lim_{n\to\infty}\int_{B_{d-1}(n)} \left\{\int_{\mathds Q_p}\chi(x_1y_1)\e^{-\sigma t|y_1|^b}\int_{\|\vec{y}_1\| \leq |y_1|} \chi(\vec{x}_1\cdot \vec{y}_1)\,{\rm d}\vec{y}_1\,{\rm d}y_1\right\}{\rm d}\vec{x}_1\end{equation} and %
\begin{equation}\label{EQ:components_marginals:I2a}I_2(x_1) = \lim_{n\to \infty}\int_{B_{d-1}(n)}\left\{\int_{\mathds Q_p}\chi(x_1y_1)\int_{\|\vec{y}_1\|>|y_1|} \chi(\vec{x}_1\cdot \vec{y}_1)\e^{-\sigma t\|\vec{y}_1\|^b}\,{\rm d}\vec{y}_1\,{\rm d}y_1\right\}{\rm d}\vec{x}_1,\end{equation} so that if both limits exist, then
\begin{equation}\label{EQ:components_marginals:I1+I2}
\rho^{(1)}(t, x_1) = I_1(x_1) + I_2(x_1).
\end{equation} %
The domain of integration of the outermost integral of both \eqref{EQ:components_marginals:I1a} and \eqref{EQ:components_marginals:I2a} is the $\mu_{d-1}$--finite measure space $B_{d-1}(n)$.
The integrand \[ \chi(x_1y_1)\int_{\|\vec{y}_1\|>|y_1|} \chi(\vec{x}_1\cdot \vec{y}_1)\e^{-\sigma t\|\vec{y}_1\|^b}\,{\rm d}\vec{y}_1 = \chi(x_1y_1)\int_{\mathds \Q_p^{d-1}} \chi(\vec{x}_1\cdot \vec{y}_1)\e^{-\sigma t\|\vec{y}_1\|^b}\mathds{1}_{\|\vec{y}_1\|>|y_1|}(\vec{y})\,{\rm d}\vec{y}_1\] is bounded with compact support in $\Q_p$, so both  \[ \chi(x_1y_1)\e^{-\sigma t|y_1|^b}\int_{\|\vec{y}_1\| \leq |y_1|} \chi(\vec{x}_1\cdot \vec{y}_1)\,{\rm d}\vec{y}_1 \quad \text{and} \quad \chi(x_1y_1)\int_{\|\vec{y}_1\|>|y_1|} \chi(\vec{x}_1\cdot \vec{y}_1)\e^{-\sigma t\|\vec{y}_1\|^b}\,{\rm d}\vec{y}_1\] are $L^1(B_{d-1}(n) \times \mathds \Q_p)$. The Fubini-Tonelli theorem implies that %
\begin{equation}\label{EQ:components_marginals:I1b}I_1(x_1) = \lim_{n\to\infty}\int_{\mathds Q_p}\left\{\int_{B_{d-1}(n)}\chi(x_1y_1)\e^{-\sigma t|y_1|^b}\int_{\|\vec{y}_1\| \leq |y_1|} \chi(\vec{x}_1\cdot \vec{y}_1)\,{\rm d}\vec{y}_1\,{\rm d}\vec{x}_1\right\}{\rm d}y_1\end{equation} and %
\begin{equation}\label{EQ:components_marginals:I2b}I_2(x_1) = \lim_{n\to \infty}\int_{\mathds Q_p}\left\{\int_{B_{d-1}(n)}\chi(x_1y_1)\int_{\|\vec{y}_1\|>|y_1|} \chi(\vec{x}_1\cdot \vec{y}_1)\e^{-\sigma t\|\vec{y}_1\|^b}\,{\rm d}\vec{y}_1\,{\rm d}\vec{x}_1\right\}{\rm d}y_1.\end{equation}

Lemma~\ref{lem:NormProd:CharInt} and \eqref{EQ:components_marginals:I1b} together imply that %
\begin{equation*}%
I_1(x_1) = \lim_{n\to\infty}\int_{\mathds Q_p} \chi(x_1y_1)\e^{-\sigma t|y_1|^b}p^{\log_p|y_1| +\min(-\log_p|y_1|, n)}{\rm d}y_1.
\end{equation*} 
For any natural number $M$, decompose the integral over $\mathds Q_p$ into a sum of integrals over $B(-M)$ and $B(-M)^c$ to obtain the equalities
\begin{align}\label{EQ:components_marginals:I1c}%
I_1(x_1) &= \lim_{n\to\infty}\left\{\int_{B(-M)^c} \chi(x_1y_1)\e^{-\sigma t|y_1|^b}p^{\log_p|y_1| +\min(-\log_p|y_1|, n)}{\rm d}y_1\right.\notag\\&\left.\qquad\qquad\qquad\qquad+ \int_{B(-M)} \chi(x_1y_1)\e^{-\sigma t|y_1|^b}p^{\log_p|y_1| +\min(-\log_p|y_1|, n)}{\rm d}y_1\right\}\notag\\%
 &= \int_{B(-M)^c} \chi(x_1y_1)\e^{-t|y_1|^b}\,{\rm d}y_1\notag\\&\qquad\qquad\qquad\qquad+ \lim_{n\to\infty}\int_{B(-M)} \chi(x_1y_1)\e^{-\sigma t|y_1|^b}p^{\log_p|y_1| +\min(-\log_p|y_1|, n)}{\rm d}y_1\notag\\
 &= \int_{{\mathds Q}_p} \chi(x_1y_1)\e^{-\sigma t|y_1|^b}\,{\rm d}y_1 + E_M(x_1),
 \end{align} 
where %
\begin{align*}
E_M(x_1) &= - \int_{B(-M)} \chi(x_1y_1)\e^{-\sigma t|y_1|^b}\,{\rm d}y_1\\& \qquad\qquad\qquad\qquad+ \lim_{n\to\infty}\int_{B(-M)} \chi(x_1y_1)\e^{-\sigma t|y_1|^b}p^{\log_p|y_1| +\min(-\log_p|y_1|, n)}{\rm d}y_1.
\end{align*}
The inequality \[|E_M(x_1)| \leq 2p^{-M}\] and \eqref{EQ:components_marginals:I1c} together imply that %
\begin{equation}\label{EQ:components_marginals:I1d}
I_1(x_1) = \int_{{\mathds Q}_p} \chi(x_1y_1)\e^{-\sigma t|y_1|^b}\,{\rm d}y_1.
\end{equation}

Take $F$ to be the function that is given for each pair $(y_1, \vec{y}_1)$ in ${\mathds Q}_p\times{\mathds Q}_p^{d-1}$ by %
\begin{equation}\label{EQ:components_marginals:F}
F(y_1, \vec{y}_1) = \e^{-\sigma t\|\vec{y}_1\|^b}{\mathds 1}_{B_{d-1}(\log_p|y_1|)^c}(\vec{y}_1)
\end{equation} 
and for each $y_1$, take $\tilde{F}(y_1, \cdot)$ to be the Fourier transform, taken over ${\mathds Q}_p^{d-1}$, of $F(y_1, \cdot)$.  Use \eqref{EQ:components_marginals:I2b} and \eqref{EQ:components_marginals:F} to rewrite $I_2(x_1)$ in terms of $F$ and obtain the equalities %
\begin{align*}
I_2(x_1) & = \lim_{n\to \infty}\int_{B_{d-1}(n)}\left\{\int_{\mathds Q_p}\chi(x_1y_1)\int_{\mathds{Q}_p^{d-1}} \chi(\vec{x}_1\cdot \vec{y}_1)F(y_1, \vec{y}_1)\,{\rm d}\vec{y}_1\,{\rm d}y_1\right\}{\rm d}\vec{x}_1\\ %
 & = \lim_{n\to \infty}\int_{B_{d-1}(n)}\left\{\int_{\mathds Q_p}\chi(x_1y_1)\tilde{F}(y_1, \vec{x}_1)\,{\rm d}y_1\right\}{\rm d}\vec{x}_1.%
 \end{align*}%
 Since $B_{d-1}(n)$ has finite measure, the Fubini-Tonelli theorem together with the equality \[F(y_1, \vec{0}) = 0\] implies that
 \begin{align}\label{EQ:components_marginals:I2PreSplit}
I_2(x_1)  & = \lim_{n\to \infty}\int_{\mathds Q_p}\left\{\int_{B_{d-1}(n)}\chi(x_1y_1)\tilde{F}(y_1, \vec{x}_1){\rm d}\vec{x}_1\right\}\,{\rm d}y_1\notag\\%
 & = \lim_{n\to \infty}\int_{\mathds Q_p}\left\{\int_{{\mathds Q}_p^{d-1}}\chi(x_1y_1)\tilde{F}(y_1, \vec{x}_1){\rm d}\vec{x}_1 - \int_{B_{d-1}(n)^c}\chi(x_1y_1)\tilde{F}(y_1, \vec{x}_1){\rm d}\vec{x}_1\right\}\,{\rm d}y_1\notag\\
 & = \lim_{n\to \infty}\int_{\mathds Q_p}\left\{F(y_1, \vec{0}) - \int_{B_{d-1}(n)^c}\chi(x_1y_1)\tilde{F}(y_1, \vec{x}_1){\rm d}\vec{x}_1\right\}\,{\rm d}y_1\notag\\
 & = -\lim_{n\to \infty}\int_{\mathds Q_p}\left\{\int_{B_{d-1}(n)^c}\chi(x_1y_1)\tilde{F}(y_1, \vec{x}_1){\rm d}\vec{x}_1\right\}\,{\rm d}y_1.
\end{align}
For any natural number $M$, decompose the integral over $\mathds Q_p$ in \eqref{EQ:components_marginals:I2PreSplit} into a sum of integrals over $B(-M)$ and $B(-M)^c$ to obtain the equality
 \begin{align}\label{EQ:components_marginals:I2Split}
I_2(x_1)  & =  -\lim_{n\to \infty}\left\{\int_{B(-M)}\left\{\int_{B_{d-1}(n)^c}\chi(x_1y_1)\tilde{F}(y_1, \vec{x}_1){\rm d}\vec{x}_1\right\}\,{\rm d}y_1 \right.\notag\\&\left.\qquad\qquad\qquad\qquad+ \int_{B(-M)^c}\left\{\int_{B_{d-1}(n)^c}\chi(x_1y_1)\tilde{F}(y_1, \vec{x}_1){\rm d}\vec{x}_1\right\}\,{\rm d}y_1\right\}.
\end{align}
For any non-zero $y_1$ in ${\mathds Q}_p$, $F(y_1, \cdot)$ is locally constant with a radius of local constancy equal to $|y_1|$.  The function $\tilde{F}(y_1, \cdot)$ is therefore supported on the ball of radius $\frac{1}{|y_1|}$ that is centered at the origin, and so the second summand in \eqref{EQ:components_marginals:I2Split} is the integral of the zero function for any $n$ that is greater than $M$.  Since the innermost integral of the first term is bounded, there is a constant $C$ so that for any natural number $M$, \[|I_2(x_1)| \leq \tfrac{C}{p^M},\] and so 
\begin{equation}\label{EQ:components_marginals:I20}
I_2(x_1) = 0.
\end{equation}

Equations~\eqref{EQ:components_marginals:I1+I2}, \eqref{EQ:components_marginals:I1d}, and \eqref{EQ:components_marginals:I20} together imply that %
\begin{equation*} 
\rho^{(1)}(t, x_1) = \int_{{\mathds Q}_p} \chi(\vec{x}_1\cdot \vec{y}_1)\e^{-\sigma t|y_1|^b}\,{\rm d}y_1.
\end{equation*}

\end{proof}

\subsection{Probabilities for the Conditioned Components}\label{sec:components:sub:probabilities}  For each $i$ in $\{1, \dots, d\}$, denote by $\vec{X}_{t, i}$ the $d-1$ tuple \[\vec{X}_{t, i} = \Big(X_t^{(1)}, \dots, X_t^{(i-1)}, X_t^{(i+1)}, \dots, X_t^{(d)}\Big), \quad \text{where}\quad \vec{X}_t =  \Big(X_t^{(1)}, \dots, X_t^{(d)}\Big).\] The components of the max-norm process fail to be independent because the spatial dependence of the law for the max-norm process involves only the component with the largest $p$-adic absolute value.  For any integers $r$ and $R$, for any $i$ in $\{1, \dots, d\}$, and for any $a$ in $B(R)$, if $r$ is less than $R$, then \[P^d\!\left(X_t^{(1)} \in B(r, a) \Big\vert \vec{X}_{t,1}\in S_{d-1}(R)\right) \leq p^{r-R}<\tfrac{1}{p}\] because this conditional probability is independent of $a$. However, as long as $t$ is small enough, \[P^d\!\left(X_t^{(1)} \in B(r, 0)\right) > \tfrac{1}{p},\] and so $\vec{X}_t$ does not have independent components as long as $t$ is small enough.  Lemma~\ref{lemma:components:conditional_calculation} provides an explicit calculation of certain conditional probabilities that lead not only to a proof that the components of $\vec{X}_t$ are dependent for any positive $t$, but also to an explicit description of certain local (small time) behaviors of the conditioned component processes.

\begin{lemma}\label{lemma:components:conditional_calculation}
For any integers $r$ and $R$, for any $i$ in $\{1, \dots, d\}$, and for any $a$ in $B(R)$, if $r$ is less than or equal to $R$, then
\begin{equation*}
P^d\!\left(X_t^{(i)} \in B(r, a) \Big\vert \vec{X}_{t,i}\in S_{d-1}(R)\right) = \frac{p^r\sum_{j\leq -R} \left(\e^{-\sigma tp^{jb}} - \e^{-\sigma tp^{(j+1)b}}\right)p^{dj}}{\sum_{j \leq -R} \left(\e^{-\sigma tp^{jb}} - \e^{-\sigma tp^{(j+1)b}}\right)p^{(d-1)j}}.
\end{equation*}
\end{lemma}

\begin{proof}%
Without loss in generality, take $i$ to be equal to $1$.  For any $\vec{x}$ in $B(r, a)\times S_{d-1}(R)$ and for any $\vec{x}$ in $B(R, a)\times S_{d-1}(R)$, \[\|\vec{x}\| = p^R,\] and so if $-j$ is less than $R$, then 
\begin{equation}\label{EQ:components:conditional_calculation:DomAZero}
\int_{B(r, a)\times S_{d-1}(R)} {\mathds 1}_{B_d(-j)}(\vec{x})\,{\rm d}\vec{x} = \int_{B(R, a)\times S_{d-1}(R)} {\mathds 1}_{B_d(-j)}(\vec{x})\,{\rm d}\vec{x} = 0.
\end{equation}
For any $\vec{x}$ in $B(R, a)^c\times S_{d-1}(R)$, \[\|\vec{x}\| > p^R,\] and so if $-j$ is less than or equal to $R$, then 
\begin{equation}\label{EQ:components:conditional_calculation:DomBZero}
\int_{B(R, a)^c\times S_{d-1}(R)} {\mathds 1}_{B_d(-j)}(\vec{x})\,{\rm d}\vec{x} = 0.
\end{equation}
The Fubini-Tonelli theorem and \eqref{EQ:NormProd:Ball_Sphere_measure_d} together imply that if \[j \leq -R,\] then
\begin{align}\label{EQ:components:conditional_calculation:DomA}
\int_{B(r, a)\times S_{d-1}(R)} {\mathds 1}_{B_d(-j)}(\vec{x})\,{\rm d}\vec{x} & = \int_{B(r, a)} {\mathds 1}_{B(-j)}(x_1)\,{\rm d}x_1 \int_{S_{d-1}(R)} {\mathds 1}_{B_{d-1}(-j)}(\vec{x}_1)\,{\rm d}\vec{x}_1\notag\\&=p^rp^{(d-1)R}\left(1-\tfrac{1}{p^{d-1}}\right),
\end{align}
and
\begin{align}\label{EQ:components:conditional_calculation:DomB}
\int_{B(R,a)\times S_{d-1}(R)} {\mathds 1}_{B_d(-j)}(\vec{x})\,{\rm d}\vec{x} & = \int_{B(R, a)} {\mathds 1}_{B(-j)}(x_1)\,{\rm d}x_1 \int_{S_{d-1}(R)} {\mathds 1}_{B_{d-1}(-j)}(\vec{x}_1)\,{\rm d}\vec{x}_1\notag\\&=p^{dR}\left(1-\tfrac{1}{p^{d-1}}\right),
\end{align}
and if \[j <-R,\] then
\begin{align}\label{EQ:components:conditional_calculation:DomC}
\int_{B(R, a)^c\times S_{d-1}(R)} {\mathds 1}_{B_d(-j)}(\vec{x})\,{\rm d}\vec{x} & = \int_{B(R, a)^c} {\mathds 1}_{B(-j)}(x_1)\,{\rm d}x_1 \int_{S_{d-1}(R)} {\mathds 1}_{B_{d-1}(-j)}(\vec{x}_1)\,{\rm d}\vec{x}_1\notag\\&=\left(p^{-j} - p^R\right)p^{(d-1)R}\left(1-\tfrac{1}{p^{d-1}}\right).
\end{align}

The equality \eqref{pdf:NormProcess:d} implies that for any positive real $t$ and any Borel subset $U$ of $\mathds Q_p^d$,
\begin{equation}\label{EQ:components:probU}
P^d(X_t\in U) = \sum_{j\in \mathds Z} \left(\e^{-tp^{jb}} - \e^{-\sigma tp^{j+1}b}\right)p^{dj}\int_U {\mathds 1}_{B_d(-j)}(\vec{x})\,{\rm d}\vec{x}.
\end{equation}
Together with \eqref{EQ:components:conditional_calculation:DomAZero}, \eqref{EQ:components:conditional_calculation:DomBZero}, \eqref{EQ:components:conditional_calculation:DomA}, \eqref{EQ:components:conditional_calculation:DomB}, and \eqref{EQ:components:conditional_calculation:DomC},  \eqref{EQ:components:probU} implies that %
\begin{equation}\label{EQ:components:CondNum}
P^d\!\left(\vec{X}_t \in B(r, a)\times S_{d-1}(R)\right) = p^rp^{(d-1)R}\left(1-\tfrac{1}{p^{d-1}}\right)\sum_{j\leq -R} \left(\e^{-\sigma tp^{jb}} - \e^{-\sigma tp^{(j+1)b}}\right)p^{dj},
\end{equation}
\begin{equation}\label{EQ:components:CondDenomA}
P^d\!\left(\vec{X}_t \in B(R, a)\times S_{d-1}(R)\right) = p^{dR}\left(1-\tfrac{1}{p^{d-1}}\right)\sum_{j\leq -R} \left(\e^{-\sigma tp^{jb}} - \e^{-\sigma tp^{(j+1)b}}\right)p^{dj},
\end{equation}
and
\begin{align}\label{EQ:components:CondDenomB}
&P^d\!\left(\vec{X}_t \in B(R, a)^c\times S_{d-1}(R)\right)  \notag\\& \hspace{1in} = p^{(d-1)R}\left(1-\tfrac{1}{p^{d-1}}\right)\sum_{j < -R} \left(\e^{-\sigma tp^{jb}} - \e^{-\sigma tp^{(j+1)b}}\right)p^{dj}\left(p^{-j} - p^R\right).
\end{align}

The equality 
\begin{align*}P\!\left(\vec{X}_{t,1}\in S_{d-1}(R)\right) &= P\!\left(\big(X_t^{(1)} \in B(R)\big)\cap\big(\vec{X}_{t,1}\in S_{d-1}(R)\big)\right) \\&\qquad\qquad+ P\!\left(\big(X_t^{(1)} \in B(R)^c\big)\cap\big(\vec{X}_{t,1}\in S_{d-1}(R)\big)\right)
\end{align*}
implies that %
\begin{align*}
&P^d\!\left(X_t^{(1)} \in B(r, a) \Big\vert \vec{X}_{t,1}\in S_{d-1}(R)\right)  \\&\qquad =  \frac{P^d\!\left(\big(X_t^{(1)} \in B(r, a)\big) \cap \big(\vec{X}_{t,1}\in S_{d-1}(R)\big)\right)}{P^d\!\left(\big(X_t^{(1)} \in B(R)\big)\cap\big(\vec{X}_{t,1}\in S_{d-1}(R)\big)\right) + P^d\!\left(\big(X_t^{(1)} \in B(R)^c\big)\cap\big(\vec{X}_{t,1}\in S_{d-1}(R)\big)\right)},
\end{align*}
and so \eqref{EQ:components:CondNum}, \eqref{EQ:components:CondDenomA}, and \eqref{EQ:components:CondDenomB} together imply that 
\begin{align*}
&P^d\!\left(X_t^{(1)} \in B(r, a) \Big\vert \vec{X}_{t,1}\in S_{d-1}(R)\right)\\&\quad = \frac{p^r\sum_{j\leq -R} \left(\e^{-\sigma tp^{jb}} - \e^{-\sigma tp^{(j+1)b}}\right)p^{dj}}{p^{R}\sum_{j\leq -R} \left(\e^{-\sigma tp^{jb}} - \e^{-\sigma tp^{(j+1)b}}\right)p^{dj} + \sum_{j < -R} \left(\e^{-\sigma tp^{jb}} - \e^{-\sigma tp^{(j+1)b}}\right)p^{dj}\left(p^{-j} - p^R\right)}%
\\&\quad = \frac{p^r\sum_{j\leq -R} \left(\e^{-\sigma tp^{jb}} - \e^{-\sigma tp^{(j+1)b}}\right)p^{dj}}{p^{R}\left(\e^{-\sigma tp^{-Rb}} - \e^{-\sigma tp^{(-R+1)b}}\right)p^{-dR} + \sum_{j < -R} \left(\e^{-\sigma tp^{jb}} - \e^{-\sigma tp^{(j+1)b}}\right)p^{(d-1)j}}%
\\&\quad = \frac{p^r\sum_{j\leq -R} \left(\e^{-\sigma tp^{jb}} - \e^{-\sigma tp^{(j+1)b}}\right)p^{dj}}{\sum_{j \leq -R} \left(\e^{-\sigma tp^{jb}} - \e^{-\sigma tp^{(j+1)b}}\right)p^{(d-1)j}}.
\end{align*}

\end{proof}

\subsection{Component Dependency}\label{sec:components:sub:dependency}  The formula that Lemma~\ref{lemma:components:conditional_calculation} provides for the conditional probabilities is rather complicated.  Lemma~\ref{lemma:components:conditional_calculation_Oest} gives a description of the local behavior of these conditional probabilities that is especially useful for understanding the effect of conditioning on the component processes.  Denote by $\Gamma(p, b, d)$ the quantity 
\begin{equation}\label{EQ:components:FDef} 
\Gamma(p, b, d) = \frac{p^{b+d}-p}{p^{b+d+1} - p}.
\end{equation}
Use the standard ``little oh'' and ``big oh'' Landau notation to simplify the statements and proofs of the statements below.

\begin{lemma}\label{lemma:components:conditional_calculation_Oest}
For any integers $r$ and $R$, for any $i$ in $\{1, \dots, d\}$, and for any $a$ in $B(R)$, if $r$ is less than or equal to $R$, then
\begin{equation*}
P^d\left(X_t^{(i)} \in B(r, a) \Big\vert \vec{X}_{t,i}\in S_{d-1}(R)\right) = \left(\Gamma(p, b, d)p^{-R} + {\rm o}(t)\right)p^r.
\end{equation*}
\end{lemma}

\begin{proof}
Take $G(R, d, \cdot)$ to be the function that is given for any $t$ in $[0, 1]$ by %
\begin{equation}\label{EQ:components:Gdef}
G(R, d, t) = \sum_{j\leq -R} \big(\e^{-\sigma tp^{jb}} - \e^{-\sigma tp^{(j+1)b}}\big)p^{dj}.
\end{equation}
Differentiate $G(R, d, \cdot)$ to obtain for each $t$ in $[0, \infty)$ the equality
\begin{align}\label{EQ:components:NumDenom}
G^\prime(R, d, 0) = \sum_{j\leq -R} \big(\sigma p^{(j+1)b} - \sigma p^{jb}\big)p^{dj} = \frac{\sigma\big(p^b-1\big)p^{(b+d)(1-R)}}{p^{(b+d)} - 1}.
\end{align}
The twice continuous differentiability of $G(R, d, \cdot)$ implies that 
\begin{equation}
G(R, d, t) = G^\prime(R, d, 0)t + {\rm O}(t^2)\quad {\rm and}\quad G(R, d-1, t) = G^\prime(R, d-1, 0)t + {\rm O}(t^2),
\end{equation}
and so if $t$ is positive, then
\begin{align}\label{EQ:components:CondbigOestimate} 
\frac{G(R,d,t)}{G(R,d-1,t)} &= \frac{G^\prime(R, d, 0)}{G^\prime(R, d-1, 0)} + {\rm o}(t)\notag\\& = p^{-R}\frac{p^{b+d}-p}{p^{b+d+1} - p} + {\rm o}(t).
\end{align}
Use \eqref{EQ:components:CondbigOestimate} to rewrite the righthand side of the equality that is given by Lemma~\ref{lemma:components:conditional_calculation} and obtain the desired description of the local behavior of the conditional probabilities.
\end{proof}

For any positive real number $t$ and any integer $R$, take $U(R)_t$ to be the random variable with the following law: For any Borel subset $V$ of $\mathds Q_p$,
\begin{equation}\label{EQ:uniformBorelDef}
{\rm Prob}(U(R)_t\in V) = P^d\!\left(X_t^{(1)} \in V \Big\vert \vec{X}_{t,1}\in S_{d-1}(R)\right).
\end{equation}

\begin{proposition}\label{prop:components:epsilon}
The random variable $U(R)_t$ is asymptotically uniformly distributed in $B(R)$.  Furthermore, for any Borel subset $V$ of $B(R)$, \[\lim_{t\to 0^+} {\rm Prob}(U(R)_t\in V) = \mu(V)p^{-R}\Gamma(p, b, d).\]
\end{proposition}

\begin{proof}
For any Borel subset $V$ of $B(R)$, there is an at most countable index set $J$ and sequences $(a_j)$ in $B(R)$ and $(r_j)$ in $\mathds Z\cap(-\infty, R]$, both indexed in $J$, so that $(B(r_j, a_j))$ is a sequence of disjoint balls whose union is $V$.  Lemma~\ref{lemma:components:conditional_calculation} and \eqref{EQ:components:CondbigOestimate} together imply that for each $j$, 
\begin{align}\label{EQ:components:epsilonA}
  {\rm Prob}(U(R)_t\in B(r_j, a_j)) &= p^{r_j}\frac{G(R, d, t)}{G(R, d-1, t)}\notag\\
                                    & = p^{r_j}\left(p^{-R}F(p,b,d) + {\rm o}(t)\right).
\end{align}
The countable additivity of the conditioned measure implies that
\begin{align}\label{EQ:components:epsilonA}
{\rm Prob}(U(R)_t\in V) &= \sum_{j\in J}p^{r_j}\left(p^{-R}\Gamma(p,b,d) + {\rm o}(t)\right)\notag\\
&= \left(\Gamma(p, b, d)p^{-R} + {\rm o}(t)\right)\sum_{j\in J}\mu(B(r_j,a_j)) = \left(\Gamma(p, b, d)p^{-R} + {\rm o}(t)\right)\mu(V).
\end{align}  
Lemma~\ref{lemma:components:conditional_calculation_Oest} and \eqref{EQ:components:epsilonA} together imply that
\begin{equation}
{\rm Prob}(U(R)_t\in V) = \left(\Gamma(p, b, d)p^{-R} + {\rm o}(t)\right)\mu(V) \to \mu(V)p^{-R}\Gamma(p, b, d)
\end{equation}
as $t$ tends to $0$ from the right, and so $U(R)_t$ is asymptotically uniformly distributed in $B(R)$.
\end{proof}

The symmetry between spatial and temporal scalings for the law for a $p$-adic Brownian motion suggests the following significant extension of Proposition~\ref{prop:components:epsilon}.

\begin{theorem}\label{theorem:components:epsilon}
For any positive real numbers $t$ and $\varepsilon$, there is an integer $M$ so that for any integer $N$ and for any Borel subset $V$ of $B(N)$, 
 \[N \ge M \quad \text{implies that} \quad \left|{\rm Prob}(U(N)_t\in V) - \mu(V)\Gamma(p, b, d)p^{-N}\right| < \varepsilon.\]
\end{theorem}

\begin{proof}
For any integer $R$ and any natural number $K$, take $G(R+K, d, t)$ and $G(R+K, d-1, t)$ to be given by \eqref{EQ:components:Gdef} so that
\begin{equation}\label{EQ:theorem:components:epsilonA}
  \frac{G(R+K,d,t)}{G(R+K,d-1,t)} = \frac{\sum_{j\leq -R-K} \left(\e^{-\sigma tp^{jb}} - \e^{-\sigma tp^{(j+1)b}}\right)p^{dj}}{\sum_{j \leq -R-K} \left(\e^{-\sigma tp^{jb}} - \e^{-\sigma tp^{(j+1)b}}\right)p^{(d-1)j}}.
\end{equation}
Reindex the sums in \eqref{EQ:theorem:components:epsilonA} to obtain the equalities
\begin{align}\label{EQ:theorem:components:epsilonB}
\frac{G(R+K,d,t)}{G(R+K,d-1,t)} & = \frac{\sum_{j\leq -R} \left(\e^{-\sigma tp^{(j-K)b}} - \e^{-\sigma tp^{(j-K+1)b}}\right)p^{d(j-K)}}{\sum_{j \leq -R} \left(\e^{-\sigma tp^{(j-K)b}} - \e^{-\sigma tp^{((j-K)+1)b}}\right)p^{(d-1)(j-K)}}\notag\\
 & = \frac{\sum_{j\leq -R} \left(\e^{-\sigma (tp^{-Kb})p^{jb}} - \e^{-\sigma (tp^{-Kb})p^{(j+1)b}}\right)p^{dj}p^{-Kd}}{\sum_{j \leq -R} \left(\e^{-\sigma (tp^{-Kb})p^{jb}} - \e^{-\sigma (tp^{-Kb})p^{(j+1)b}}\right)p^{(d-1)j}p^{-K(d-1)}}\notag\\
& = p^{-K}\frac{G(R,d,tp^{-Kb})}{G(R,d-1,tp^{-Kb})}.%
\end{align}
For any positive real number $t$, Lemma~\ref{lemma:components:conditional_calculation_Oest} implies that 
\[
\frac{G(R,d,tp^{-Kb})}{G(R,d-1,tp^{-Kb})} = \left(\Gamma(p, b, d)p^{-R} + {\rm o}(tp^{-Kb})\right),
\]
and so \eqref{EQ:theorem:components:epsilonB} implies that
\begin{equation}\label{EQ:theorem:components:epsilonC}
\frac{G(R+K,d,t)}{G(R+K,d-1,t)} = p^{-K}\left(\Gamma(p, b, d)p^{-R} + {\rm o}(tp^{-Kb})\right).
\end{equation}
Lemma~\ref{lemma:components:conditional_calculation_Oest}, \eqref{EQ:uniformBorelDef}, and \eqref{EQ:theorem:components:epsilonC} together imply  that for any $r$ that is less than $R+K$ and any $a$ in $B(R+K)$,
\begin{equation*}
{\rm Prob}(U(R+K)_t\in B(r,a)) = p^{r-K}\left(\Gamma(p, b, d)p^{-R} + {\rm o}(tp^{-Kb})\right).
\end{equation*}
Follow the proof of Proposition~\ref{prop:components:epsilon} to generalize to the case where $V$ is a Borel set that is not a ball and obtain for any Borel set $V$ in $B(R+K)$ the equality 
\begin{equation*}
{\rm Prob}(U(R+K)_t\in V) = \mu(V)\Gamma(p, b, d)p^{-R-K} + \mu(V)p^{-K}{\rm o}(tp^{-Kb}).
\end{equation*}
Since $V$ is a subset of $B(R+K)$, $\mu(V)$ is no greater than $p^{R+K}$, and so Proposition~\ref{prop:components:epsilon} implies that \[\left|\mu(V)p^{-K}{\rm o}(tp^{-Kb})\right| \leq p^R\left|{\rm o}(tp^{-Kb})\right| < \varepsilon\] as long as $K$ is large enough.

\end{proof}

Note that Theorem~\ref{theorem:components:epsilon} implies that for any fixed positive $t$, $U(R)_t$ is asymptotically uniformly distributed in $R$ for large values of $R$.

\begin{corollary}\label{thm:components:anyt}
For any positive real number $t$, the components of $\vec{X}_t$ are not independent. 
\end{corollary}

\begin{proof}
For any $i$ in $\{1, \dots, d\}$ and for any positive real number $t$, Theorem~\ref{thm:components:marginal_distributions} implies that 
\[
\lim_{r\to \infty} P^d\!\left(X_t^{(i)} \in B(r)\right) = \lim_{r\to \infty} P(X_t \in B(r)) = 1,
\]
and so for any positive real number $\varepsilon$ there is an integer $R$ so that for any $r$ that is greater than or equal to $R$, 
\[
P^d\!\left(X_t^{(i)} \in B(r)\right) > 1 - \varepsilon.
\]
Since $\Gamma(p,b,d)$ is less than 1, Theorem~\ref{theorem:components:epsilon} implies that there is a natural number $K$ and a real number $e$ in $(-\varepsilon, \varepsilon)$ so that
\begin{align*}
P^d\!\left(X_t^{(i)} \in B(R+K) \Big\vert \vec{X}_{t,i}\in S_{d-1}(R+1+K)\right) &= \mu(B(R+K))p^{-R-K-1}\Gamma(p, b, d)+e\\
& = \tfrac{1}{p}\Gamma(p, b, d)+e < \tfrac{1}{p} +\varepsilon.
\end{align*}
As long as $\varepsilon$ is small enough, \[P^d\!\left(X_t^{(i)} \in B(R+K) \Big\vert \vec{X}_{t,i}\in S_{d-1}(R+1+K)\right) < P^d\!\left(X_t^{(i)} \in B(R+K)\right).\]
\end{proof}


\section{First Exit Probabilities}\label{sec:exit}

Each component process of the product and max-norm process is a $p$-adic Brownian motion with diffusion constant $\sigma$ and diffusion exponent $b$. The dependency of the components of these processes impacts the first exit times from balls.  Namely, dimension has a large impact on the first exit time probabilities for the product process, but it has a rather small effect on these probabilities for the max-norm process.

\subsection{Exit Times for the Components}\label{sec:exit:sub:1D}

Take $\alpha$ to be the positive real number that is given by \[\alpha = 1- \tfrac{p^b-1}{p^{b+1}-1}\] and take $X$ to be a one dimensional $p$-adic Brownian motion with diffusion constant $\sigma$ and diffusion exponent $b$.  For sake of clarity, slightly modify the notation in \cite{Weisbart:2021} and for any positive real number $T$ denote by $\vertiii{X}_T$ the quantity \[\vertiii{X}_T=\sup_{0\leq t\leq T}|X_t|.\]  The probability that a sample path for $X$ remains in $B(R)$ until time $T$, $P\!\left(\vertiii{X}_T\leq p^R\right)$, is a \emph{survival probability} for $X$.  The complement of this probability is a \emph{first exit probability}, the probability that a sample path for $X$ has first exit from $B(R)$ before time $T$.   Since every point in $B(R)$ is the center of $B(R)$, the exit times for $X$ from $B(R)$ do not depend on starting points, and so \[P\!\left(\vertiii{X}_{T+S}\leq p^R \big\vert{X}_{S}\leq p^R \right) = P\!\left(\vertiii{X}_{T}\leq p^R\right).\]  The survival probabilities for $X$ are continuous from the right at $0$ and satisfy Cauchy's multiplicative functional equation, and so the first exit time for $X$ is an exponentially distributed random variable.  Theorem~3.1 of \cite{Weisbart:2021} determines the parameter of this exponential distribution by establishing the equality%
\begin{equation}\label{ExitProbability:Theorem_3.1}
P\!\left(\vertiii{X}_T \leq p^R\right) = {\rm e}^{-\sigma \alpha T p^{-Rb}}.
\end{equation}

\subsection{Exit Times for the Processes}\label{sec:exit:sub:exit}

Take $\alpha_d$ to be the quantity \[\alpha_d = 1- \tfrac{p^b-1}{p^{b+d}-1},\] and for any $\mathds Q_p^d$--valued stochastic process $\vec{Z}$ and any positive real number $T$, take $\vertiii{\vec{Z}}_T$ to be the quantity \[\vertiii{\vec{Z}}_T=\sup_{0\leq t\leq T}\|\vec{Z}_t\|.\]  

\begin{theorem}\label{theorem:exit:Prod}
For any integer $R$,
\begin{equation*}
\otimes^dP\!\left(\vertiii{\vec{Y}}_T \leq p^R\right) = {\rm e}^{-d\sigma \alpha_1 T p^{-Rb}}.
\end{equation*}
\end{theorem}

\begin{proof}
Since $\vec{Y}_t$ lies outside $B_d(R)$ if and only if at least one component of $\vec{Y}_t$ lies outside $B_d(R)$, \[\left(\vertiii{\vec{Y}}_T\leq p^R\right) = \bigcap_{i\in\{1, \dots, d\}}\left(\vertiii{Y^{(i)}}_T\leq p^R\right).\] Independence of the components of $\vec{Y}$ implies that \[\otimes^dP\!\left(\vertiii{\vec{Y}}_T\leq p^R\right) = \prod_{i\in\{1, \dots, d\}}P\!\left(\vertiii{Y^{(i)}}_T\leq p^R\right) = {\rm e}^{-d\alpha\sigma T p^{-Rb}}.\]
\end{proof}

With only minor modification, the arguments of \cite{Weisbart:2021} extend to the more general max-norm setting and determine the survival probabilities for the max-norm process $\vec{X}$.  For this reason, the proof below for Theorem~\ref{theorem:exit:max-norm} will omit certain details that the proof of Theorem~3.1 in \cite{Weisbart:2021} includes.

\begin{theorem}\label{theorem:exit:max-norm}
For any integer $R$, \[P^d\!\left(\vertiii{\vec{X}}_T\leq p^R\right) = {\rm e}^{-\sigma\alpha_d Tp^{-Rb}}.\]
\end{theorem}

\begin{proof} %

Take $U$ in \eqref{EQ:components:probU} to be $B(R)$ in order to obtain the equality
\begin{align}\label{eq:components:intoverball}
P^d\!\left(\vec{X}_t\in B_d(R)\right) & = \sum_{j\in \mathds Z} \left(\e^{-\sigma tp^{jb}} - \e^{-\sigma tp^{(j+1)b}}\right)p^{dj}\int_{B_d(R)} {\mathds 1}_{B_d(-j)}(\vec{x})\,{\rm d}\vec{x}\notag\\
& = p^{dR}\sum_{j \leq -R} \left(\e^{-\sigma tp^{jb}} - \e^{-\sigma tp^{(j+1)b}}\right)p^{dj} + \sum_{j>-R} \left(\e^{-\sigma tp^{jb}} - \e^{-\sigma tp^{(j+1)b}}\right).
\end{align}
The summands of the second sum in \eqref{eq:components:intoverball} telescope, and so
\begin{equation}\label{eq:components:intoverballSimp}
P^d\!\left(\vec{X}_t\in B_d(R)\right) = \e^{-\sigma tp^{(-R+1)b}} + p^{dR}\sum_{j\leq -R} \left(\e^{-\sigma tp^{jb}} - \e^{-\sigma tp^{(j+1)b}}\right)p^{dj}.
\end{equation}

For any natural numbers $N$ and $j$, where $j$ is less than or equal to $N$, take $t_j$ to be given by \[t_j = \tfrac{jT}{N}.\]  Since the max-norm on $\mathds Q_p^d$ satisfies the ultra-metric inequality, 
\begin{align}\label{sec3:maxtoprodequality}&P^d\!\left(\max_{t_j} \big(\|\vec{X}_{t_1}\|, \dots, \|\vec{X}_{t_n}\|\big) \leq p^R\right)\notag\\&\hspace{.5in} = P^d\!\left(\max_{t_j} \big(\|\vec{X}_{t_1}\|, \|\vec{X}_{t_1}- \vec{X}_{t_2}\|, \dots, \|\vec{X}_{t_n}- \vec{X}_{t_{n-1}}\|\big) \leq p^R\right).\end{align}  Take the random variable $X_0$ to be the zero function.  The independence of the set of increments $\big\{\vec{X}_{t_j} - \vec{X}_{t_{j-1}}\colon i\in\{1, \dots, N\}\big\}$ implies that \begin{align}\label{eq:exit:ProdInc}P^d\!\left(\max_{t_j} \big(\|\vec{X}_{t_1}\|, \dots, \|\vec{X}_{t_n}\|\big) \leq p^R\right) &= P^d\!\left(\|\vec{X}_{t_1}\| \leq p^R\cap \dots\cap \|\vec{X}_{t_n}\| \leq p^R\right)\notag\\& = \prod_{1\leq j\leq N}P^d\!\left(\|\vec{X}_{t_i}-\vec{X}_{t_{i-1}}\| \leq p^R\right)^N.\end{align}  The increments are identically distributed, and so 
\begin{equation}
 P^d\!\left(\max_{t_j} \big(\|\vec{X}_{t_1}\|, \dots, \|\vec{X}_{t_n}\|\big) \leq p^R\right)= P^d\!\left(\|\vec{X}_{\frac{T}{N}}\| \leq p^R\right)^N.
 \end{equation}
 
Take $B$ to be the twice continuously differentiable function that is defined for any $t$ in $[0,\infty)$ by \[B(t) = P^d\!\left(\|\vec{X}_t\| \leq p^R\right).\]  The twice continuous differentiability of $B$ implies that the equality \[B\big(\tfrac{T}{n}\big) = 1 + \tfrac{B^\prime(0)T}{n} + {\rm O}\big(\tfrac{1}{n^2}\big),\quad \text{hence} \quad  B\big(\tfrac{T}{n}\big)^n = {\rm e}^{TB^\prime(0)}.\] The right continuity of the sample paths of $\vec{X}$ implies that
\begin{align}\label{eq:exit:EqforSurvivalProb}P^d\!\left(\vertiii{\vec{X}}_T\leq p^R\right) &= \lim_{n\to \infty}P^d\!\left(\max_{j} \big(\|\vec{X}_{\frac{T}{n}}\|, \|\vec{X}_{\frac{2T}{n}}\|, \dots, \|\vec{X}_{\frac{jT}{n}}\|, \dots, \|\vec{X}_{T}\|\big) \leq p^R\right)\notag\\ &= \lim_{n\to \infty} B\big(\tfrac{T}{n}\big)^n = {\rm e}^{TB^\prime(0)}.\end{align}  

Rewrite \eqref{eq:components:intoverballSimp} to obtain the equality 
\begin{equation}\label{eq:exit:intoverballSimp}
B(t) = \e^{-\sigma tp^{-Rb}} + p^{dR}\sum_{j < -R} \left(\e^{-\sigma tp^{jb}} - \e^{-\sigma tp^{(j+1)b}}\right)p^{dj}.
\end{equation}
Differentiate both sides of \eqref{eq:exit:intoverballSimp} to obtain the equality%
\begin{equation}\label{eq:exit:Bprime}
B^\prime(0) = -\sigma p^{-Rb} + \sigma \sum_{j\leq -R-1}p^{d(R+j)}\Big(p^{(j+1)b} - p^{jb}\Big).
\end{equation}
Simplify the righthand side of \eqref{eq:exit:Bprime} to obtain the equality
\begin{equation}\label{eq:exit:Bprimeat0}
B^\prime(0) = -\sigma \alpha_d p^{-Rb}
\end{equation} which, together with \eqref{eq:exit:EqforSurvivalProb}, implies that \[P^d\!\left(\vertiii{\vec{X}}_T\leq p^R\right) = \e^{-\sigma \alpha_d Tp^{-Rb}}.\]
\end{proof}

Together with Theorem~\ref{theorem:exit:max-norm}, the equality 
\begin{equation*}
\lim_{d\to \infty} \alpha_d = \lim_{d\to \infty} \left(1- \tfrac{p^b-1}{p^{b+d}-1}\right) = 1
\end{equation*} %
implies that for any integer $R$, 
\begin{equation}\label{EQ:exit:normLim}
\lim_{d\to \infty} P^d\!\left(\vertiii{\vec{X}}_T\leq p^R\right) = {\rm e}^{-\sigma Tp^{-Rb}}.
\end{equation} %
In contrast, Theorem~\ref{theorem:exit:Prod} implies that %
\begin{equation}\label{EQ:exit:ProdLim}
\lim_{d\to \infty} \otimes^dP\!\left(\vertiii{\vec{Y}}_T\leq p^R\right) = \lim_{d\to \infty}{\rm e}^{-d\sigma \alpha_1Tp^{-Rb}} = 0.
\end{equation}  
This marked contrast between \eqref{EQ:exit:normLim} and \eqref{EQ:exit:ProdLim} demonstrates that exit probabilities from a fixed ball depend only to a small degree on dimension for the max-norm process, but are very sensitive to changes in dimension for the product process.


\end{document}